
%

\input amstex
\loadbold
\documentstyle{amsppt}
\NoBlackBoxes

\pagewidth{32pc}
\pageheight{44pc}
\magnification=\magstep1

\def\PP{\Bbb{P}}
\def\OO{\Cal{O}}
\def\Alb{\operatorname{Alb}}
\def\Pic{\operatorname{Pic}}
\def\Proj{\operatorname{Proj}}
\def\Sym{\operatorname{Sym}}
\def\rank{\operatorname{rank}}
\def\rest#1#2{\left.{#1}\right\vert_{{#2}}}
\def\QED{{\unskip\nobreak\hfil\penalty50\quad\null\nobreak\hfil
{$\square$}\parfillskip0pt\finalhyphendemerits0\par\medskip}}

\topmatter
\title
Remarks on rational points of varieties
whose cotangent bundles are generated by global sections
\endtitle
\rightheadtext{}
\author Atsushi Moriwaki \endauthor
\leftheadtext{}
\address
Department of Mathematics, Faculty of Science,
Kyoto University, Kyoto, 606-01, Japan
\endaddress
\curraddr
Department of Mathematics, University of California,
Los Angeles, 405 Hilgard Avenue, Los Angeles, California 90024, USA
\endcurraddr
\email moriwaki\@math.ucla.edu \endemail
\date September, 1994 \enddate
\abstract
In this short note, we will gives several remarks on rational points
of varieties whose cotangent bundles are generated by global sections.
For example, we will show that
if the sheaf of differentials $\Omega^1_{X/k}$ of
a projective variety $X$ over a number field $k$ is ample and
generated by global sections, then the set of $k$-rational
points of $X$ is finite.
\endabstract
\endtopmatter

\document

\head
\S0. Introduction
\endhead

In \cite{Fa1} and \cite{Fa2},
G. Faltings proved the following theorem.

\proclaim{Theorem A}
Let $A$ be an abelian variety over a number field $k$ and
$X$ a subvariety of $A$. Then, there are a finite number of
translated abelian subvarieties $B_1, \ldots, B_n$ over $k$ such that
$B_i \subset X$ and the closure $\overline{X(k)}$ of $X(k)$ in $X$
is contained in $\bigcup_i B_i$.
\endproclaim

In this short note, as applications of the above Faltings' theorem,
we will give several remarks on rational points
of varieties whose cotangent bundles are generated by global sections.
The first remark is the following theorem, which is
a slight generalization of Theorem~A.

\proclaim{Theorem B}
Let $X$ be a projective variety over a number field $k$,
$A$ an abelian variety over $k$, and $\alpha : X \to A$ a morphism over $k$.
If $\alpha^*(\Omega^1_{A/k}) \to \Omega^1_{X/k}$ is surjective,
then every irreducible component of $\overline{X(k)}$ is
geometrically irreducible and isomorphic to an abelian variety.
\endproclaim

As corollary of Theorem~B, we have the following.

\proclaim{Corollary C}
Let $X$ be a connected smooth projective variety over a number field $k$.
If the sheaf of differentials $\Omega^1_{X/k}$ is generated by global sections,
then every irreducible component of $\overline{X(k)}$ is
geometrically irreducible and isomorphic to an abelian variety.
\endproclaim

Thus, if we assume further that $\Omega^1_{X/k}$ is ample in Corollary~C,
then $X$ contains no abelian variety. Hence,
the set of $k$-rational points of $X$ is finite.
This is a partial answer of the following conjecture
under the additional assumption ``$\Omega^1_{X/k}$ is generated by global
sections''.

\definition{Conjecture D} (cf. \cite{La})
Let $X$ be a smooth projective variety over a number field $k$.
If the cotangent bundle $\Omega^1_{X/k}$ is ample, then $X(k)$ is finite.
\enddefinition

It is true even on a singular variety, namely,

\proclaim{Theorem E}
Let $X$ be a projective variety over a number field $k$.
If the sheaf of differentials $\Omega^1_{X/k}$ of $X$ over $k$ is ample and
generated by global sections, then the set of $k$-rational
points of $X$ is finite, where
ampleness of $\Omega^1_{X/k}$ means that
the tautological line bundle on
$\Proj( \bigoplus_{d \geq 0} \Sym^d(\Omega^1_{X/k}))$ is ample.
\endproclaim

Finally, I would like to thank Prof. Lazarsfeld for his kind suggestions.

\head
\S1. Non-denseness
\endhead

In this section, we will give a powerful and simple non-denseness theorem.

First of all, we fix notation.
Let $X$ be a connected smooth projective variety over
a field $k$ of characteristic zero.
Assume $X(k) \not= \emptyset$ and fix $x_0 \in X(k)$.
Let $\Alb_{X/k}$ be the dual abelian variety of the Picard variety
$\Pic^0_{X/k}$
of $X$. Let $Q$ be the universal line bundle on
$X \times \Pic^0_{X/k}$ with
$\rest{Q}{\{x_0\}\times\Pic^0_{X/k}} = \OO_{\Pic^0_{X/k}}$.
$Q$ gives rise to an $X$-valued point of $\Alb_{X/k}$, that is,
a morphism $\alpha : X \to \Alb_X$.
Actually, $\alpha$ is given by
$\alpha(x) = \rest{Q}{\{x\}\times\Pic^0_{X/k}}$ for $x \in X$.
In particular, $\alpha(x_0) = 0$.
The abelian variety $\Alb_{X/k}$ is called the Albanese variety of $X$ over $k$
and
$\alpha$ is called the Albanese map. Moreover,
the dimension of $\Alb_{X/k}$, denoted by $q(X)$,
is called the irregularity of $X$.
Further, we denote by $d(X)$ the dimension of
$\alpha(X)$, i.e. $d(X) = \dim \alpha(X)$.
In general, $d(X) \leq q(X)$ and $d(X) \leq \dim X$.
It is well known that $q(X) = \dim_k H^0(X, \Omega^1_{X/k}) = \dim_k H^1(X,
\OO_X)$.

Let $Y$ be a geometrically irreducible projective variety over $k$
and $\mu : X \to Y$ a desingularization of $Y$.
It is easy to see that $q(X)$ and $d(X)$ does not depend on
the choice of $X$. Thus, we can define $q(Y)$ and $d(Y)$
by $q(X)$ and $d(X)$ respectively.

As a consequence of Theorem~A, we have the following non-denseness theorem.

\proclaim{Theorem 1.1}
Let $X$ be a geometrically irreducible projective variety
over a number field $k$.
If $d(X) < q(X)$,
then the closure of $X(k)$ in $X$
is a proper closed subset.
\endproclaim

\demo{Proof}
Considering a desingularization of $X$, we may assume that
$X$ is smooth over $k$. Moreover, we may assume $X(k) \not= \emptyset$.
Let $\alpha : X \to \Alb_{X/k}$ be the Albanese map.
Since $q(X)  > d(X)$, $\alpha$ is not surjective.
If $X(k)$ is dense in $X$,
then so is $\alpha(X)(k)$ in $\alpha(X)$. Thus, by Theorem~A,
$\alpha(X)$ is an abelian subvariety of $\Alb_{X/k}$.
This is a contradiction because
$\alpha(X)(\bar{k})$ generates $\Alb_{X/k}(\bar{k})$.
\QED
\enddemo

\head
\S2. Lemmas
\endhead

In this section, we will prepare three lemmas for the proof of Theorem~E.

Let $X$ be a projective scheme and $\Cal{F}$
a coherent $\OO_X$-module. We denote
$\Proj( \bigoplus_{d \geq 0} \Sym^d(\Cal{F}) )$ by $\PP(\Cal{F})$.
Let $\pi_{F} : \PP(\Cal{F}) \to X$ be the natural morphism and
$\OO_{\Cal{F}}(1)$ the tautological line bundle on $\PP(\Cal{F})$.
We say $\Cal{F}$ is {\it ample} if $\OO_{\Cal{F}}(1)$ is ample.

\proclaim{Lemma 2.1}
Let $X$ be a projective scheme and $\Cal{F}$
a coherent $\OO_X$-module.
\roster
\item "(1)" Let $\Cal{F} \to \Cal{G}$ be a surjective homomorphism
of $\OO_X$-modules. If $\Cal{F}$ is ample, then so is $\Cal{G}$.

\item "(2)" Let $f: Y \to X$ be a finite morphism of projective schemes.
If $\Cal{F}$ is ample, then so is $f^*(\Cal{F})$.
\endroster
\endproclaim

\demo{Proof}
(1) Since $\Cal{F} \to \Cal{G}$ is surjective,
$\PP(\Cal{G})$ is a subscheme of $\PP(\Cal{F})$ and
$\OO_{\Cal{G}}(1) = \rest{\OO_{\Cal{F}}(1)}{\PP(\Cal{G})}$.
Thus, $\OO_{\Cal{G}}(1)$ is ample.

(2) Let $f' : \PP(f^*(\Cal{F})) \to \PP(\Cal{F})$ be the induced morphism.
Since $\OO_{f^*(\Cal{F})}(1) = {f'}^*(\OO_{\Cal{F}}(1))$ and
$f'$ is finite, $\OO_{f^*(\Cal{F})}(1)$ is ample.
\QED
\enddemo

\proclaim{Lemma 2.2}
Let $X$ be a scheme over a field $k$ and
$\Cal{F}$ a coherent $\OO_X$-module.
\roster
\item "(1)" Let $\alpha : \Cal{F} \to \Cal{G}$ be a surjective homomorphism
of $\OO_X$-modules.
If $\Cal{F}$ is generated by global sections, then so is $\Cal{G}$.

\item "(2)" Let $f: Y \to X$ be a morphism of schemes over $k$.
If $\Cal{F}$ is generated by global sections, then so is $f^*(\Cal{F})$.

\item "(3)" If $\Cal{F}$ is generated by global sections,
then so is $\OO_{\Cal{F}}(1)$.
\endroster
\endproclaim

\demo{Proof}
(1) Let us consider the following commutative diagram:
$$
\CD
H^0(X, \Cal{F}) \otimes_k \OO_X @>>> H^0(X, \Cal{G}) \otimes_k \OO_X \\
@VVV  @VVV \\
\Cal{F} @>>> \Cal{G}
\endCD
$$
Since $H^0(X, \Cal{F}) \otimes_k \OO_X \to \Cal{F}$ and $\Cal{F} \to \Cal{G}$
are surjective, so is $H^0(X, \Cal{G}) \otimes_k \OO_X \to \Cal{G}$.

(2) Since $H^0(X, \Cal{F}) \otimes_k \OO_X \to \Cal{F}$ is surjective,
we have a surjection $H^0(X, \Cal{F}) \otimes_k \OO_Y \to f^*(\Cal{F})$.
Here we consider the following diagram:
$$
\CD
H^0(X, \Cal{F}) \otimes_k \OO_Y @>>> f^*(\Cal{F}) \\
@VVV @| \\
H^0(Y, f^*(\Cal{F})) \otimes_k \OO_Y @>>> f^*(\Cal{F})
\endCD
$$
Thus, $H^0(Y, f^*(\Cal{F})) \otimes_k \OO_Y @>>> f^*(\Cal{F})$ is surjective.

(3) By (2), $\pi_{\Cal{F}}^*(\Cal{F})$ is generated by global sections.
On the other hand, there is the natural surjective homomorphism
$\pi_{\Cal{F}}^*(\Cal{F}) \to \OO_{\Cal{F}}(1)$.
Thus, by (1), $\OO_{\Cal{F}}(1)$ is generated by global sections.
\QED
\enddemo

\proclaim{Lemma 2.3}
Let $Y$ be a geometrically irreducible projective variety
over a field $k$ of characteristic zero.
If $\dim Y \geq 1$ and $\Omega^1_{Y/k}$ is ample and generated by global
sections,
then $q(Y) \geq 2\dim Y$.
\endproclaim

\demo{Proof}
Let $\mu : Y' \to Y$ be a desingularization of $Y$ such that
$\mu^*(\Omega^1_{Y/k})/\mu^*(\Omega^1_{Y/k})_{tor}$ is locally free.
We set $P = \PP(\mu^*(\Omega^1_{Y/k})/\mu^*(\Omega^1_{Y/k})_{tor})$ and
$L = \OO_{\mu^*(\Omega^1_{Y/k})/\mu^*(\Omega^1_{Y/k})_{tor}}(1)$.
By Lemma~2.2, $\mu^*(\Omega^1_{Y/k})/\mu^*(\Omega^1_{Y/k})_{tor}$
is generated by global sections.
Thus, so is $L$ by Lemma~2.2.
Hence we have a morphism $\phi_{|L|} : P \to \PP^N$
with $\phi_{|L|}^*(\OO_{\PP^N}(1)) = L$, where $N = \dim_k H^0(P, L) - 1$.
Let $\nu$ be the composition of maps
$P \hookrightarrow \PP(\mu^*(\Omega^1_{Y/k})) \to \PP(\Omega^1_{Y/k})$.
Then, $L = \nu^*(\OO_{\Omega^1_{Y/k}}(1))$.
Since $\nu$ gives a birational morphism
from $X$ to the image $\nu(X)$ and $\OO_{\Omega^1_{Y/k}}(1)$ is ample,
$L$ is big.
It follows that $\phi_{|L|}$ is generically finite. Therefore,
$$
   \dim P = \dim \phi_{|L|}(P) \leq \dim \PP^N = \dim_k H^0(P, L) - 1,
$$
which implies
$$
\dim_k H^0(Y', \mu^*(\Omega^1_{Y/k})/\mu^*(\Omega^1_{Y/k})_{tor})  =
\dim_k H^0(P, L) \geq \dim P + 1 = 2 \dim Y.
$$
On the other hand, the natural homomorphism
$\mu^*(\Omega^1_{Y/k}) \to \Omega^1_{Y'/k}$ induces an injection
$\mu^*(\Omega^1_{Y/k})/\mu^*(\Omega^1_{Y/k})_{tor} \to \Omega^1_{Y'/k}$.
Hence,
$$
\dim_k H^0(Y', \Omega^1_{Y'/k}) \geq
\dim_k H^0(Y', \mu^*(\Omega^1_{Y/k})/\mu^*(\Omega^1_{Y/k})_{tor}).
$$
Therefore, $\dim_k H^0(Y', \Omega^1_{Y'/k}) \geq 2 \dim Y$, which says
$q(Y) \geq 2 \dim Y$.
\QED
\enddemo

\head
\S3. Proofs of Theorem~B, Corollary~C and Theorem~E
\endhead

\subhead \S3.1. Proof of Theorem~B \endsubhead
Let $Y$ be an irreducible component of $\overline{X(k)}$.
Since $Y(k)$ is dense in $Y$, $Y$ is geometrically irreducible.
For, let $Y_{\bar{k}} = \Gamma_1 \cup \cdots \cup \Gamma_r$ be
the irreducible decomposition of $Y_{\bar{k}}$.
Then, it is easy to see $Y(k) \subset \bigcap_i \Gamma_i$.
Thus, if $r \geq 2$, we have a contradiction.

We set $B = \alpha(Y)$.
Since $Y(k)$ is dense in $Y$, so is $B(k)$ in $B$.
Thus, by Faltings' theorem (Theorem~A), $B$ is a translated abelian
subvariety of $A$. Let us consider the following diagram:
$$
\CD
\rest{\alpha^*(\Omega^1_{A/k})}{Y} @>{\beta}>> \rest{\Omega^1_{X/k}}{Y} \\
@V{\delta}VV @VV{\epsilon}V \\
\left(\rest{\alpha}{Y}\right)^*(\Omega^1_{B/k}) @>>{\gamma}> \Omega^1_{Y/k} \\
\endCD
$$
Since $\beta$, $\delta$ and $\epsilon$ are surjective,
so is $\gamma$. On the other hand, $\rank \Omega^1_{B/k} \leq
\rank \Omega^1_{Y/k}$ and $\Omega^1_{B/k}$ is locally free.
It follows that $\gamma$ gives an isomorphism between
$\left(\rest{\alpha}{Y}\right)^*(\Omega^1_{B/k})$ and
$\Omega^1_{Y/k}$. Thus, $Y$ is smooth over $k$ and
$\rest{\alpha}{Y}$ is etale.
Therefore, by a theorem due to S. Lang (cf. \cite{Mu, Chapter IV, 18}),
$Y$ is an abelian variety.
\QED

\subhead \S3.2. Proof of Corollary~C \endsubhead
Let us consider the Albanese map $\alpha : X \to \Alb_{X/k}$.
Since $H^0(X, \Omega^1_{X/k}) \otimes \OO_X \to \Omega^1_{X/k}$
is surjective and $H^0(X, \Omega^1_{X/k}) \otimes \OO_{\Alb_{X/k}}
\simeq \Omega^1_{\Alb_{X/k}}$, $\alpha^*(\Omega_{\Alb_{X/k}}) \to
\Omega^1_{X/k}$ is surjective. Therefore, we can apply Theorem~B
to get our corollary.
\QED

\subhead \S3.3. Proof of Theorem~E \endsubhead
Assume that $\overline{X(k)}$ has an irreducible component
$Y$ with $\dim Y \geq 1$.
In the same way as in the proof of Theorem~B,
$Y$ is geometrically irreducible.
By (2) of Lemma~2.1 and (2) of Lemma~2.2, $\rest{\Omega^1_{X/k}}{Y}$
is ample and generated by global sections. Here there is the natural surjection
$\rest{\Omega^1_{X/k}}{Y} \to \Omega^1_{Y/k}$. Thus, by virtue of (1) of
Lemma~2.1 and
(1) of Lemma~2.2, $\Omega^1_{Y/k}$ is ample and generated by global sections.
Therefore, by Lemma~2.3, $q(Y) \geq 2 \dim Y$.
Hence, by Theorem~1.1,
$Y(k)$ is not dense in $Y$. This is a contradiction.
\QED

\subhead \S3.4. Remark \endsubhead
Over a function field, a positive answer of Conjecture~D was obtained by
Noguchi
(cf. \cite{No} and \cite{Mo}). In this case, the set of rational points
is however not finite in general. They are concentrated on a proper closed
subset.


\enddocument